\documentclass[11pt]{article}
\usepackage{amssymb,amsfonts,amsmath,latexsym,epsf,tikz,url}
\newtheorem{theorem}{Theorem}[section]

\newtheorem{lemma}[theorem]{Lemma}

\newtheorem{definition}[theorem]{Definition}

\newcommand{\proof}{\noindent{\bf Proof.\ }}
\newcommand{\qed}{\hfill $\square$\medskip}

\begin{document}

\title{Some new results on the  total domination polynomial of a graph}

\author{Saeid Alikhani$^{}$\footnote{Corresponding author}  and  Nasrin Jafari }

\date{\today}

\maketitle

\begin{center}

   Department of Mathematics, Yazd University, 89195-741, Yazd, Iran\\
{\tt alikhani@yazd.ac.ir}\\

\end{center}


\begin{abstract}
Let $G = (V, E)$ be a simple graph of order $n$. The total dominating set of $G$ is a subset $D$ of $V$ that every vertex of $V$ is adjacent to some vertices of $D$. The total domination number of $G$ is equal to minimum cardinality of  total dominating set in $G$ and is denoted by $\gamma_t(G)$. The total domination polynomial of $G$ is the polynomial $D_t(G,x)=\sum_{i=\gamma_t(G)}^n d_t(G,i)x^i$, where $d_t(G,i)$ is the number of total dominating sets of $G$ of size $i$. A root of $D_t(G,x)$ is called a total domination root of $G$. An irrelevant edge of $D_t(G,x)$ is an edge $e \in E$,  such that $D_t(G, x) = D_t(G\setminus e, x)$. In this paper, we characterize edges possessing this property. Also  we obtain some results for the number of total dominating sets of a regular graph. Finally, we study graphs with exactly two total domination roots $\{-3,0\}$, $\{-2,0\}$ and $\{-1,0\}$.  
\end{abstract}

\section{Introduction}

Let $G = (V, E)$ be a simple graph. The order of $G$ is the number of vertices of $G$. For any vertex $ v \in V$, the open neighborhood of $v$ is the set $N(v)=\{ u \in V | uv \in E\}$ and the closed neighborhood is the set $N[v]=N (v) \cup \{v\}$.
For a set $S\subset V$, the open neighborhood of $S$ is the set $N(S)=\bigcup_{v\in S }N(v)$ and the closed neighborhood of $S$ is the set $N[S]=N (S) \cup S$. The set $D\subset V$ is a total dominating set if every vertex of $V$ is adjacent to some vertices of $D$, or equivalently, $N(D)=V$. The total domination  number $\gamma_t(G)$ is the minimum cardinality of a total dominating set in $G$. A total dominating set with cardinality $\gamma_t(G)$ is called a $\gamma_t$-set. An $i$-subset of $V$ is a subset of $V$ of cardinality $i$. Let $\mathcal{D}_t(G, i)$ be the family of total dominating sets of $G$ which are $i$-subsets and let $d_t(G,i)=|\mathcal{D}_t(G, i)|$. The polynomial $D_t(G; x)=\sum_{i=1}^n d_t(G,i)x^i$ is defined as total domination polynomial of $G$.
A root of $D_t(G, x)$ is called a total  domination root of $G$. We denote the set of distinct total domination roots by $Z(D_t(G,x))$.

The corona of two graphs $G_1$ and $G_2$, as defined by Frucht and Harary in \cite{harary}, is the graph $ G_1\circ G_2$ formed from one copy of $G_1$ and $|V(G_1)|$ copies of $G_2$, where the $i$-th vertex of $G_1$ is adjacent to every vertex in the $i$-th copy of $G_2$. The corona $G\circ K_1$, in particular, is the graph constructed from a copy of $G$, where for each vertex $v \in V(G)$, a new vertex $v'$ and a pendant edge $vv'$ are added.

Recurrence relations of graph polynomials have received considerable attention in the
literature. It is well-known that the independence
polynomial and matching polynomial of a graph satisfies a linear recurrence relation with respect to two vertex elimination operations, the deletion of a vertex and the deletion of vertex's closed neighborhood. Other
 graph polynomials in the literature satisfy similar recurrence relations with
respect to vertex and edge elimination operations \cite{Kotek}. 
In contrast, it is significantly harder to find recurrence relations for the domination
polynomial and the total domination polynomial.  The easiest recurrence relation is to remove an edge and to compute the total domination polynomial of the graph arising instead of the one for the original graph. Indeed, for the total domination polynomial of a graph there might be such irrelevant edges, that can be deleted without changing the value of the total domination polynomial at all.  
 An irrelevant edge is an edge $e \in E$ of $G$, such that $D_t(G, x) = D_t(G\setminus e, x)$.

The roots of graph polynomials reflect some important information about the structure of graphs. There are many papers on the location of the roots of graph polynomials such as chromatic polynomial, matching polynomial, independence polynomial, characteristic polynomial and domination polynomial. We refer the reader to \cite{Oboudi} and its references for more information in roots of graph polynomials.   In \cite{nasrin} we have shown that all roots of $D_t(G, x)$ lie in the circle with center $(-1, 0)$ and the
radius $\sqrt[\delta]{2^n-1}$, where $\delta$ is the minimum degree of $G$. Also we proved that for a graph $G$ of order $n$, if $\delta\geq \frac{2n}{3}$,
then every integer root of $D_t(G, x)$ lies in the set $\{-3,-2,-1,0\}$.

\noindent As usual we denote the complete graph, path and cycle of order $n $ by $K_n$, $P_n$ and $C_n$, respectively. Also $S_n$ is the star graph with $n$ vertices. A leaf (end-vertex) of a graph is a vertex of degree one, while a support vertex is a vertex adjacent to a leaf.

\medskip

In the next section,  we characterize irrelevant edges for the total domination polynomial. We consider regular graphs in Section 3 and study their total domination polynomials.   Finally we study graphs with exactly two total domination roots $\{-3,0\}$, $\{-2,0\}$ and $\{-1,0\}$ in Section 4. 


\section{Irrelevant edges}

The easiest recurrence relation for total domination polynomial of a graph is to remove an edge and to compute the total domination polynomial of the graph arising instead of the one for the original graph. Indeed, for the total domination polynomial of a graph there might be such irrelevant edges, that can be deleted without changing the value of the total domination polynomial at all. In this section, we study these  edges.

\begin{definition}
Let $G = (V,E)$ be a graph. A vertex $v\in V$  is total domination-covered, if each total dominating set of $G\setminus v$ is a total dominating set of $G$ and the total dominating sets of $G\setminus v$ are exactly those total dominating sets of $G$ not including the vertex $v$.
\end{definition} 

The proofs of Theorems \ref{V1}  and \ref{E1} follow the proofs in \cite{Kotek} with some minor changes:

\begin{theorem}\label{V1}
Let $G = (V,E)$ be a graph. A vertex $v \in  V$ is total domination-covered if there is a vertex $u \in N[v]$ such that $N[u] \subseteq N[v]$.
\end{theorem}

\proof
 To dominate $u$, the vertex $u$ or a vertex adjacent to $u$ must be in each total dominating set of $G\setminus v$. Since  every vertex adjacent to $u$ is also adjacent to $v$ in $G$, so we have the result.\qed

\begin{theorem}\label{E1}
Let $G=(V,E)$ be a graph. If $e=\{u, v\}\in  E$ is an irrelevant edge in $G$, then $u$ and $v$ are total domination-covered in $G \setminus  e$.
\end{theorem}

\proof
By contradiction, suppose that  at least one vertex (say $u$) is not total domination-covered in  $G \setminus  e$.  Then there exist a total dominating set of $G\setminus e$ 
which  is not a total dominating set of $G$, and this implies that $D_t(G,x)\neq D_t(G\setminus e,x)$ which is a contradiction. \qed

\begin{figure}
\begin{center}
\includegraphics[width=3.2cm,height=2.7cm]{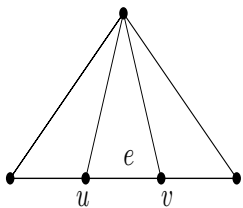}
\caption{The edge $e$ is not irrelevant but $u$ and $v$ are total domination covered in $G\setminus e$. }
\label{G}
\end{center}
\end{figure}

Note that the converse of Theorem \ref{E1} is not  true. As an example  for  the graph in Figure \ref{G}, the vertices  $u$ and $v$ are total domination-covered,  while $e$ is not irrelevant edge. Because $D_t(G,x)=x^5+5x^4+8x^3+5x^2$ but the total domination polynomial of $G\setminus e$ is $ x^5+5x^4+6x^3+4x^2$.

We need the following theorem to obtain more results. 
\begin{theorem}{\rm\cite{Dod}} \label{D}
If  $G=(V,E)$ is  a graph and   $e=\{u,v\}\in E $ with  $N[v]=N[u]$, then 
             $D_t(G, x)=D_t(G\setminus e, x)+x^2D_t(G\setminus N[u],x)$.
      \end{theorem}

 By Theorem \ref{D}, we have the following result.
 
 \begin{theorem}
 Let $G=(V, E)$ be a graph and   $e=\{u,v\}\in E $ with  $N[v]=N[u]$. If there exists  a support vertex $w\in N[u]$, then the edge  $e$ is an irrelevant edge.
 \end{theorem}
 
 \proof
 By theorem \ref{D}, we have
  \begin{center}
 $D_t(G, x)=D_t(G\setminus e, x)+x^2D_t(G\setminus N[u],x)$.
 \end{center}
 Note that the graph $G\setminus N[u]$ has at least an isolated vertex, therefore $D_t(G\setminus N[u],x)=0$ and we have the result.\qed
 
\begin{theorem}\label{irrelevant}
Let $G$ be a graph and $e=\{u,v\}$ is an edge of $G$. If the vertices $u$ and $v$ are adjacent to the support vertices, then $e$ is an irrelevant edge.
\end{theorem}

\proof
Suppose that the vertices  $u$, $v$ are adjacent to the  support vertices   $w$ and $z$, respectively. Then, every total dominating set of $G$ include support vertices  $w$ and $z$. So the vertices on edge $e$, under any total dominating set of $G$ are dominated and adjacency between them is ineffective. Therefore $D_t(G,x)=D_t(G\setminus e,x)$ and $e$ is an irrelevant edge.\qed

  Let to compute the total domination polynomial of  a  family of graphs which has shown in Figure \ref{FNK} using the irrelevant edges. 
   An $(n,k)$-firecracker $F(n,k)$ is a graph obtained by the concatenation of $n$, $k$-stars $S_k$ by linking one leaf from each. See figure \ref{FNK}. The following easy theorem gives the total domination number of this kind of  graphs: 
  
\begin{figure}
\begin{center}
\includegraphics[width=9.75cm,height=2.6cm]{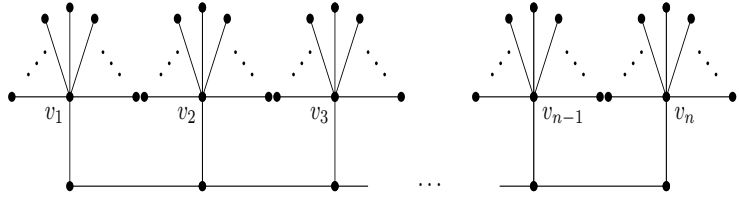}
\caption{The graph $F(n,k)$.}
\label{FNK}
\end{center}
\end{figure}

\begin{theorem}
For every natural numbers $n$ and $k$, we have $\gamma_t(F(n,k))=2n$.
\end{theorem}
 
 \proof
 Let $D$ be a minimum total dominating set of $F(n,k)$. Then  $\{v_1, v_2, \ldots, v_n\}\subseteq D$ and for every $v_i $, the set $D$ contains exactly one vertex that is adjacent to $v_i$. So $\gamma_t(F(n,k))=2n$.
 \qed
 
\begin{theorem}
For every natural numbers $n$ and $k\geq 3$, 
\begin{center}
  $D_t(F(n,k),x)=(x(x+1)^{(k-1)}-x)^n$.
 \end{center}
\end{theorem}

\proof
By Theorem \ref{irrelevant}, every edge that linking $k$-stars together is an irrelevant edge. Therefore $D_t(F(n,k),x)=(D_t(S_k,x))^n$ and we have the  result.\qed

\begin{figure}
\begin{center}
\includegraphics[width=9.25cm, height=2.25cm]{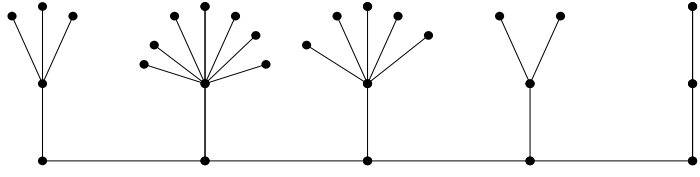}
\caption{The graph $F(5,9,7,4,3)$}
\label{F}
\end{center}
\end{figure}

Now, we generalize the definition of firecracker graphs. An $(k_1,k_2,\ldots,k_n)$-firecracker $F(k_1,\ldots, k_n)$ is a a graph obtained by the concatenation of $k_i$-stars $S_{k_i}$ by linking one leaf from each (see Figure \ref{F}). If $k_i\geq 3$  ($1\leq i\leq n$), then  every edge that linking $k_i$-stars's together, is an irrelevant edge. Therefore $D_t(F(k_1,\ldots, k_n))=\prod \limits_{i=1}^n (x(x+1)^{(k_i-1)}-x)$. 

\medskip

Here, we are interested to examine the effect on the total domination polynomial of a graph when we remove a vertex. Recall that a vertex $v\in G$ is called essential vertex, if $D_t(G\setminus v,x)=0$ (\cite{Dod}). 

\begin{lemma}
Let $G=(V,E)$ be a graph. The vertex $v\in V(G)$ is an essential vertex if and only if $v$ is a support vertex of $G$.
\end{lemma}
\proof  Since  $D_t(G,x)=0$ if and only if $G$ has an isolated vertex, so we have the result.\qed


\begin{theorem}{\rm\cite{Dod}}\label{Dod}
If  $G=(V,E)$ is  a graph and $v\in V(G)$, then 
\begin{center}
$D_t(G,x)=D_t(G\setminus v,x)+D_t(G\odot v,x)-D_t(G\circledcirc v,x).$
\end{center}
where $G\odot v$ is the graph obtained from $G$ by removing all edges between vertices of $N(v)$ and $G\circledcirc v=G\odot v\setminus v$.
\end{theorem}

\begin{lemma}
Let $G=(V,E)$ be a graph and  $v$ is a support  vertex of $G$. Then $D_t(G,x)=D_t(G\odot v,x)$.
\end{lemma}
\proof
Since the vertex $v$ is a support vertex,  so $D_t(G\setminus v,x)=D_t(G\circledcirc v,x)=0$. Therefore,  by Theorem \ref{Dod}, we have the result.\qed

\section{Total domination polynomial of regular graphs}

In this section, we study  some coefficients of the total domination polynomial of regular graphs and then compute the total domination polynomial of cubic graphs  of order $10$. 

We denote the family of all total dominating  sets of $G$ with cardinality $i$ and contain a vertex $v$ by $\mathcal{D}_t^v(G,i)$ and $d_t^v(G,i)=|\mathcal{D}_t^v(G,i)|$. Two graphs $G$ and $H$ are said to be total dominating equivalent or simply $\mathcal{D}_t$-equivalent, if $D_t(G,x)=D_t(H,x)$ and written $G\thicksim H$. It is evident that the relation $\thicksim$ of $\mathcal{D}_t$-equivalent is an equivalence relation on the family $\mathcal{G}$ of graphs, and thus $\mathcal{G}$ is partitioned into equivalence classes, called the $\mathcal{D}_t$-equivalence classes. Given $G\in \mathcal{G}$, let 
\begin{center}
$[G]=\{H\in \mathcal{G} ~: H\thicksim G\}.$
\end{center}

If $[G]=\{G\}$, then $G$ is said to be total dominating unique or simply $\mathcal{D}_t$-unique. 
In this section, similar to \cite{Cubic}, we determine the $\mathcal{D}_t$-equivalence classes for cubic graphs of order $10$. The proofs of Theorems \ref{dtv}  and \ref{regular}
  follow the proofs in \cite{Cubic} with some minor changes:

\begin{figure}[!h]
	\hglue2.5cm
	\includegraphics[width=11cm,height=5cm]{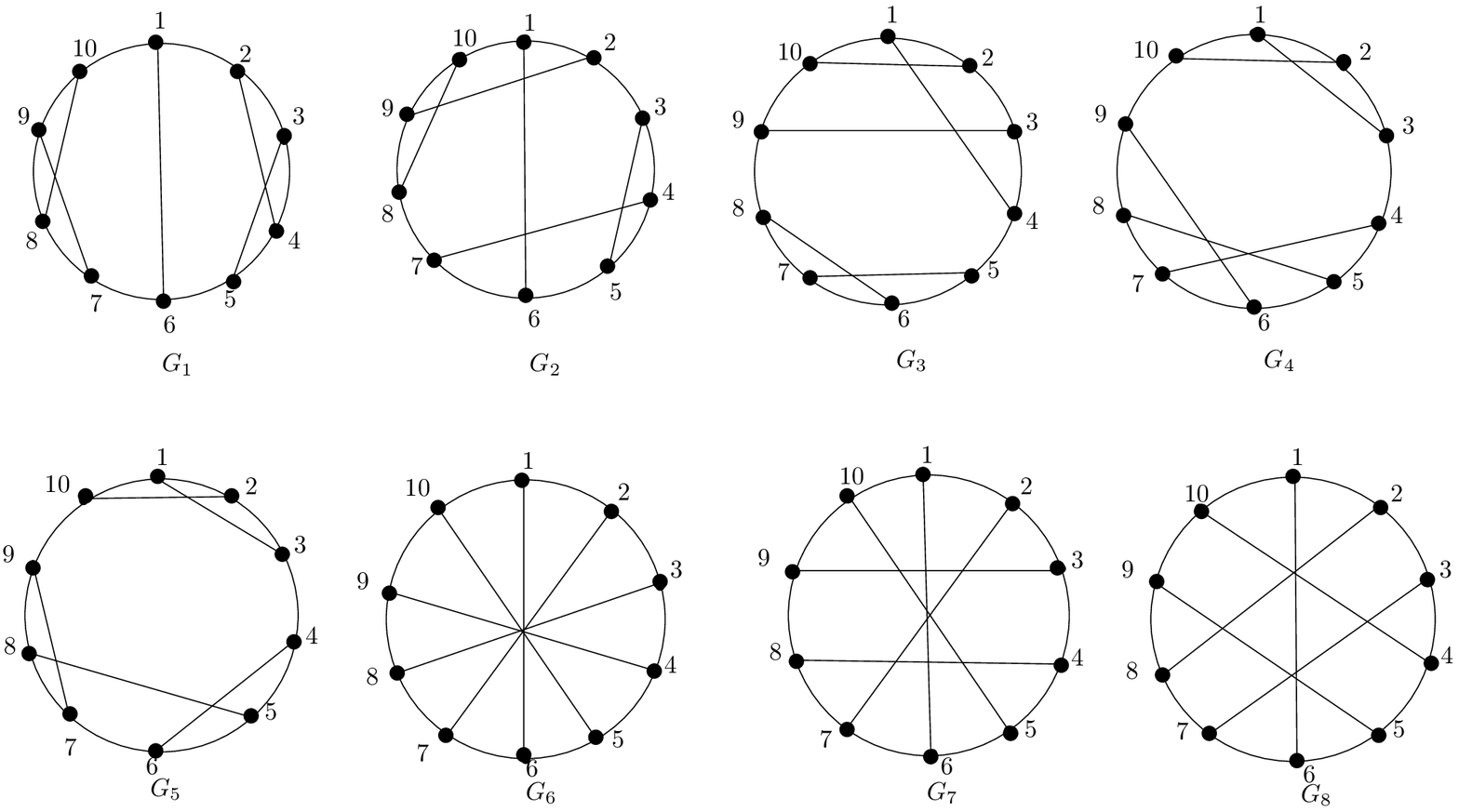}
	\vglue5pt
	\hglue2.5cm
	\includegraphics[width=11cm,height=5cm]{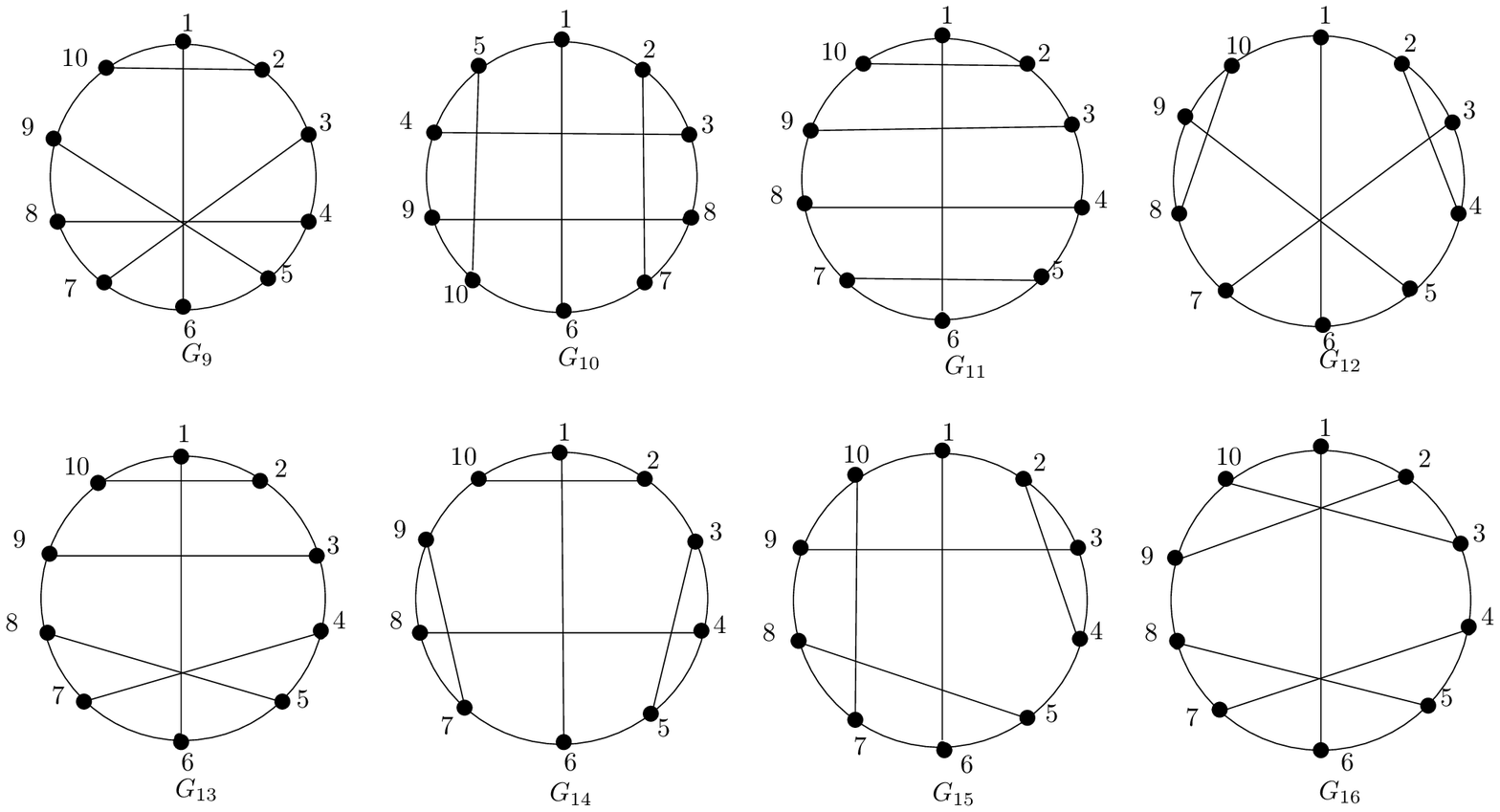}
	\hglue2.5cm
	\includegraphics[width=10.7cm,height=5cm]{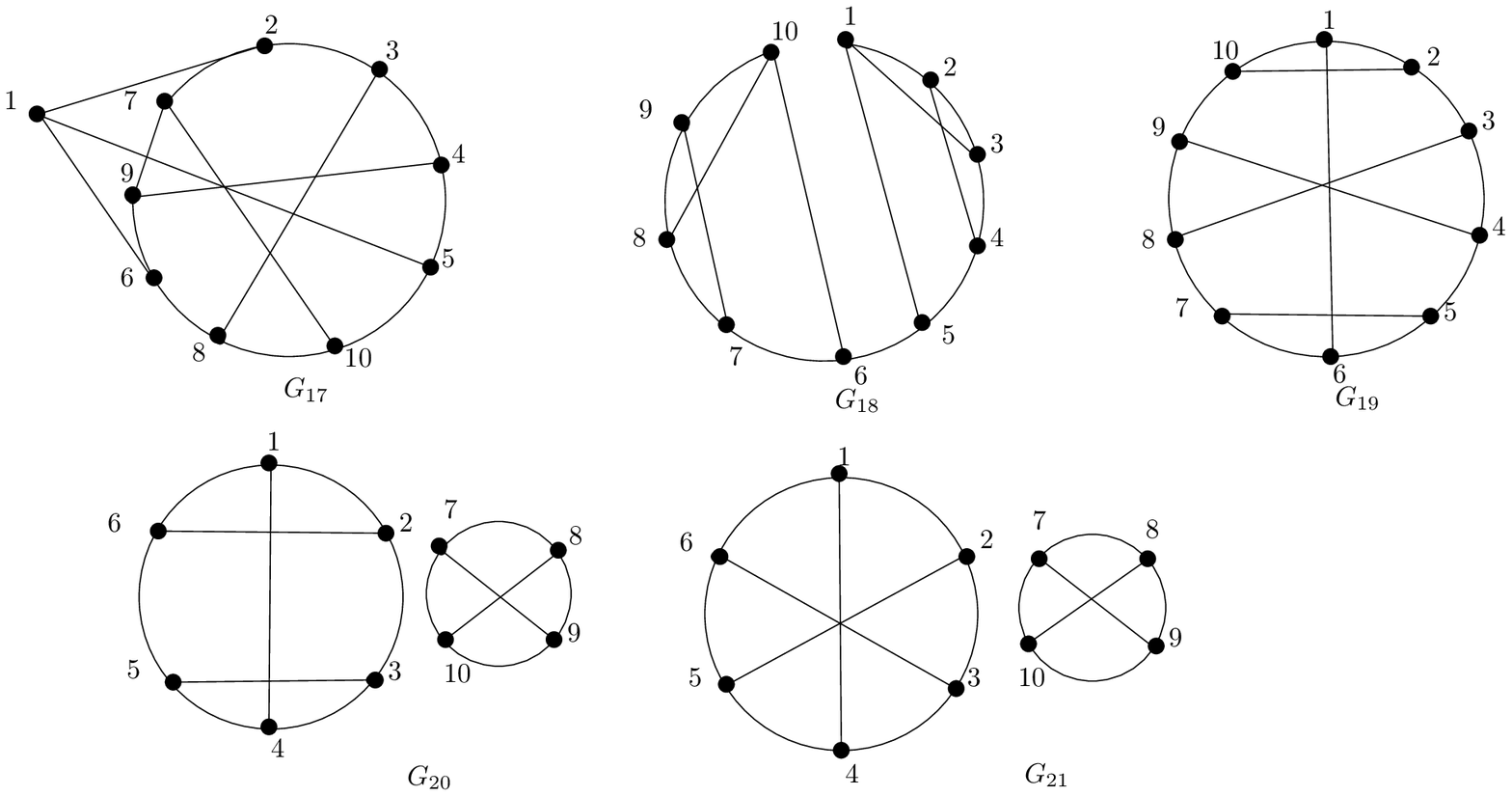}
	\vglue-10pt \caption{\label{figure2} Cubic graphs of order $10$.}
\end{figure}

\begin{lemma}\label{dtv}
Let $G=(V,E)$ be a vertex transitive graph of order $n$ and $v\in V$. For any $1\leq i\leq n$, we have $d_t(G,i)=\frac{n}{i}d_t^v(G,i)$.
\end{lemma}

\proof
If $D$ is a total dominating set of vertex transitive graph, $G$, with size $i$ and $\theta \in Aut(G)$, then $\theta(D)$ is also a total dominating set of $G$ with size $i$. Also, because $G$ is a vertex transitive graph, so for every vertices $v$ and $u$, $d_t^v(G,i)=d_t^u(G,i)$. If $D$ is a total dominating set of size $i$, then there are exactly $i$ vertices $v_{j_1}, v_{j_2},\ldots, v_{j_i}$ such that $D$ counted in $d_t^{v_{j_k}}(G,i)$, for any $1\leq k\leq i$. Therefore $d_t(G,i)=\frac{n}{i}d_t^v(G,i)$ and the proof is complete.\qed

\begin{lemma}\label{regular}
Let $G=(V,E)$ be $k$-regular graph of order $n$. Then $d_t(G,i)=\binom{n}{i}$ for all $i>n-k$.
\end{lemma}

\proof
Let $G$ be a $k$-regular graph of order $n$ with vertex set $V$. Each vertex of $G$ is adjacent to $k$ vertices. Let $S$ be a  $(k-1)$-subsets of $V$, and  $V'=V\setminus S$. Then  $V'$ is a total dominating set for $G$ of size $n-k+1$ and the number of total dominating set of size $n-k+1$ for $G$ is equal to number of ways of choosing $k-1$ vertex of $V$. Therefore $d_t(G,n-k+1)=\binom{n}{k-1}=\binom{n}{n-k+1}$. Similarly for each $2\leq i\leq k-1$, we have $d_t(G,n-i)=\binom{n}{n-i}$. So for any $i>n-k$, $d_t(G,i)=\binom{n}{i}$ and the proof is complete.\qed

In the study of the total domination polynomial of regular graphs, it is natural to ask about the total domination polynomial of Petersen graph and its $\mathcal{D}_t$-equivalence class. To answer to this question, we consider exactly $21$ cubic graphs of order $10$ given in Figure \ref{figure2} (see \cite{Cubic}). There are just two non-connected cubic graphs of order $10$.  Note that the graph $G_{17}$ is the Petersen graph. 
 The following theorem gives the total domination polynomial of the Petersen graph. 

\begin{theorem}
The total domination polynomial of Petersen graph $P$ is $$D_t(P,x)=x^{10}+10x^9+45x^8+110x^7+140x^6+72x^5+10x^4.$$
\end{theorem}

\proof
We have $\gamma_t(P)=4$. Since the Petersen graph $P$ is a $3$-regular graph of order $10$, by Lemma \ref{regular}, we have $d_t(P,i)=\binom{10}{i}$, for $i=8,9,10$. On the other hand, since $P$ is a vertex transitive graph, we calculate $d_t(P,i)$, for $i=4,5,6$ using Lemma \ref{dtv}. So we have the result. \qed
 

Using Maple we computed the total domination polynomial of  cubic graphs of order $10$. As some consequences we stat the following results for graphs in Figure \ref{figure2}. 
\begin{theorem}
	\begin{enumerate} 
		\item[(i)] The Petersen graph $P$ is not $\mathcal{D}_t$-unique. More precisely, the three graphs $G_{12}$, $G_{14}$ and $P\cong G_{17}$ are $\mathcal{D}_t$-equivalent.
\item[(ii)] 
The three  graphs $G_1$, $G_8$ and $G_9$ are $\mathcal{D}_t$-equivalent.
\item[(iii)]
The graphs $G_2$, $G_3$, $G_4$, $G_5$, $G_6$, $G_7$, $G_{10}$, $G_{11}$, $G_{13}$, $G_{15}$, $G_{16}$, $G_{18}$, $G_{20}$, $G_{21}$  are $\mathcal{D}_t$-unique.
\end{enumerate} 
\end{theorem}

\section{On the graphs with exactly two total domination roots }

Graphs whose certain  polynomials have few roots can sometimes give interesting  information about the structure of  graph. The characterization of graphs with few distinct roots of characteristic polynomials (i.e., graphs with few distinct eigenvalues) have been the subject of many researchers \cite{bri,dam1,dam2,dam3}. Also the first authors 
has studied graphs with few domination roots in \cite{few}.  In \cite{nasrin} we have shown that all roots of $D_t(G, x)$ lie in the circle with center $(-1, 0)$ and the
radius $\sqrt[\delta]{2^n-1}$, where $\delta$ is the minimum degree of $G$. Also  we proved that for a graph $G$ of order $n$, if $\delta\geq \frac{2n}{3}$,
then every integer root of $D_t(G, x)$ lies in the set $\{-3,-2,-1,0\}$.  Motivated by these integer roots, and a conjecture in \cite{nasrin} which states that for every  integer root
$r$ of $D_t(G,x)$,  $r\in\{-3,-2,-1,0\}$, 
  we study graphs with exactly two total domination roots $\{-1,0\}$, $\{-2,0\}$ and $\{-3,0\}$, in this section.

\subsection{Graphs with exactly two total domination roots $\{-1,0\}$}
In this subsection, first we state and prove the following theorem to present a necessary condition  for a graph to have exactly two total domination roots $-1$ and $0$.

\begin{theorem}\label{R}
If  $G=(V,E)$ is  a graph of order $n$ with $r$ support  vertices, then $d_t(G,n-1)=n-r$.
\end{theorem}

\proof
Let $A\subseteq V$ be the set of all support vertices of $G$. For every vertex $v\in V\setminus A$, the set $V\setminus \{v\}$ is a total dominating set of  $G$. So $d_t(G,n-1)=n-r$.\qed

\begin{theorem}
Let $G=(V,E)$ be a graph of order $n$. If $Z(D_t(G,x))=\{-1,0\}$, then $G$ has at least two  support vertices.
\end{theorem}

 \proof Let $G$ be a graph of order $n$ and $D_t(G,x)=x^a(x+1)^b$, such that $a+b=n$ and $a=\gamma_t(G)\geq 2$. By Theorem \ref{R}, $d_t(G,n-1)=n-r$, where $r$ is the number of support vertices. So $n-r=b$ and $a=r$. Therefore we have result.\qed

\subsection{Graphs with exactly two total domination roots $\{-2,0\}$}

\begin{figure}[h]
	\begin{center}
		\includegraphics[width=14cm, height=3.5cm]{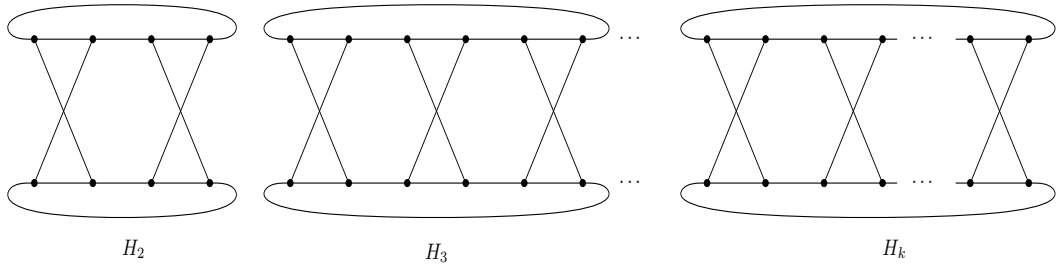}
		\caption{The graphs $\mathcal{H}$.}
		\label{mathcalH}
	\end{center}
\end{figure}

In this subsection, we present a necessary condition for  graphs  two total domination roots $\{-2,0\}$.  
We recall that  a vertex cut  of a graph $G$ is a subset $V'$ of $V(G)$ such that $G-V'$ is not connected and a  $k$-vertex cut of $G$ is a vertex cut of $k$ vertices.
The connectivity, $\kappa(G)$ of a connected graph $G$ (which contains no complete graph factor) is the smallest integer $k$ for which $G$ has a $k$-vertex cut.
To obtain our result, we introduce  graphs in Figure \ref{mathcalH} which denoted  by $\mathcal{H}$. 

Infinite family  $\mathcal{H}$ are connected cubic graphs. For $k \geq 2$, let $H_k$ be the graph constructed as follows. Consider two copies of the path $P_{2k}$ with respective vertex sequences $a_1b_1a_2b_2 \ldots a_kb_k$ and $c_1d_1c_2d_2 \ldots c_kd_k$. Let $A=\{a_1,a_2,\ldots ,a_k\}$, $B=\{b_1,b_2,\ldots ,b_k\}$, $C=\{c_1,c_2,\ldots ,c_k\}$ and $D=\{d_1,d_2, \ldots ,d_k\}$. For each $i\in \{1,2, \ldots ,k\}$,  join $a_i$ to $d_i$ and $b_i$ to $c_i$. To complete the construction of the graph $H_k$, join $a_1$ to $b_k$ and $c_1$ to $d_k$. We note that $H_k$ are cubic graphs of order $4k$.

\begin{figure}[h]
	\begin{center}
		\includegraphics[width=4cm, height=4cm]{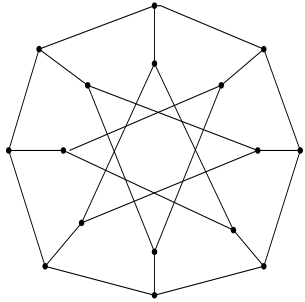}
		\caption{The generalized Petersen graph of order $16$.}
		\label{GPetersen}
	\end{center}
\end{figure}

\begin{theorem}{\rm \cite{bri}}\label{Gp}
	If $G$ is  a $3$-connected graph of order $n$, then $\gamma_t(G)\leq \frac{n}{2}$ with equality if and only if $G=K_4$ or $G\in \mathcal{H}$ or $G$ is the generalized Petersen graph $GP$ of order $16$ shown in Figure \ref{GPetersen}.
\end{theorem}

\begin{theorem}
	Let $G$ be a simple graph of order $n$. If $Z(D_t(G,x))=\{-2,0\}$, then $\kappa(G)\leq 2$.
\end{theorem}

\proof
Let $G$ be a $3$-connected graph and $D_t(G,x)=x^{\gamma_t}(x+2)^{\beta}$, where $\gamma_t+\beta=n$. So by Theorem \ref{R}, we have $d_t(G,n-1)=n=2\beta$, and so  $\beta=\frac{n}{2}$, $\gamma_t(G)=\frac{n}{2}$. By Theorem \ref{Gp},  $G=K_4$ or $G\in \mathcal{H}$ or $G$ is the generalized Petersen graph, $GP$, of order $16$. But  
$D_t(K_4,x)=x^4+4x^3+6x^2$, $D_t(GP,x)=x^8(x^4+8x^3+28x^2+48x+30)^2$, and so 
$G\neq K_4, GP$.
On the other hand, for every $k=2,\ldots,m$, $H_k\in \mathcal{H}$ is a $3$-regular graph of order $4k$ with $\gamma_t(H_k)=2k$. So by Lemma \ref{regular} we have
$d_t(H_k,4k-2)=2k(4k-1).$
By the assumption, we shall have

$$D_t(G,x)= x^{2k}(x+2)^{2k}=\sum\limits_{i=\gamma_t(G)}^{2k}\binom{2k}{i} 2^{2k-i}x^{2k+i}.$$

So $d_t(G,4k-2)=4k(2k-1)$, which  is a contradiction. Therefore $G\notin \mathcal{H}$ and we have result. \qed

\begin{figure}
	\begin{center}
		\includegraphics[width=3.7cm,height=3.2cm]{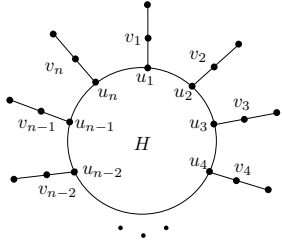}
		\caption{The graph $H(3)$.}
		\label{H(3)}
	\end{center}
\end{figure}

As an example of family of graph $G$ with $Z(D_t(G,x))=\{-2,0\}$, let $H$ be an arbitrary graph of order $n$ and consider $n$ copies of graph $P_3$.   
By definition, the graph $H(3)$ is obtained by identifying each vertex of $H$ with an end vertex of a $P_3$ (\cite{Nasrin2}). See Figure \ref{H(3)}.
By Theorem \ref{irrelevant} we compute the total domination polynomial of  $H(3)$ (see also \cite{Nasrin2}).

\begin{theorem}\label{H}
	For any graph $H$ of order $n$, we have $D_t(H(3),x)=x^{2n}(x+2)^n$.
\end{theorem}

\proof
Let $D$ be a total dominating set of $H(3)$  of size $k\geq n$ in  Figure \ref{H(3)}. Obviously $\{v_1,v_2,\ldots,v_n\}\subset D$ and for any $v_i$, the set $D$ contains exactly one vertex that is adjacent to $v_i$, so $\gamma_t(H(3))=2n$. On the other hand, for all $i,j$, $ 1\leq i\neq j\leq n$ and $e=\{u_i, u_j\}$, the vertices $u_i, u_j$ are adjacent to support vertices $v_i$, $v_j$. Therefore each  edge in $H$ is an irrelevant edge, and so we have
\begin{center}
	$D_t(H(3),x)=(D_t(P_3,x))^n=x^{2n}(x+2)^n. $
\end{center}\qed

\subsection{Graphs with exactly two total domination roots $\{-3,0\}$}
In this subsection,  we shall characterize graphs whose total domination polynomial have exactly two roots $-3$ and $0$.  To do this,  we need the following result. 

\begin{theorem}{\rm\cite{bri}}\label{3}
Let $G$ be a connected graph of order $n\geq 3$. Then $\gamma_ t (G)=\frac{2n}{3}$ if and only if $G$ is $C_3$, $C_6$ or $H(3)$ (in Figure \ref{H(3)}) for some connected graph $H$.
\end{theorem}

The next theorem classifies   all connected graphs without end-vertices,  whose total domination polynomial have just two roots $\{-3,0\}$.

\begin{theorem}
Let $G=(V,E)$ be a graph of order $n$ with $\delta(G)\geq 2$.  Then $Z(D_t(G,x))=\{-3,0\}$ if and only if $G$ is $C_3$ or $C_6$.
\end{theorem}

\proof
First  note that $D_t(C_3,x)=x^2(x+3)$ and $D_t(C_6,x)=x^4(x+3)^2$.  Let $D_t(G,x)=x^{\alpha}(x+3)^{\beta}$. Therefore $\alpha+\beta=n$ and $\alpha=\gamma_t(G)$. By Theorem \ref{R} we have $d_t(G,n-1)=n-r$ where $r$ is the number of support vertices  of $G$. Therefore $d_t(G,n-1)=n$ . On the other hand 

\begin{equation*}
 D_t(G,x)=x^{\alpha}(x+3)^{\beta}=\sum\limits_{i=0}^{\beta}\binom{\beta}{i}x^{\alpha +i}3^{\beta-i}.
\end{equation*}
Therefore  we have $d_t(G,n-1)=n=3\beta$ and so  $\gamma_t(G)=\frac{2n}{3}$. By Theorem \ref{3}, $G$ is $C_3$, $C_6$ or $H(3)$ for some connected graph $H$. Since $\delta(G)\geq 2$ and  by Theorem \ref{H}, $G$ is $C_3$ or $C_6$.\qed

\end{document}